\documentclass[11pt]{article}
\usepackage{amssymb,amsmath,amsthm}
\begin{document}
\title{Mergelyan's approximation theorem with nonvanishing polynomials and universality of zeta-functions}

\theoremstyle{plain}
\newtheorem{thm}{Theorem}
\newtheorem{lem}{Lemma}\renewcommand{\thelem}{\hskip -4 pt}
\newtheorem{cor}{Corollary}
\theoremstyle{definition} 
\newtheorem{example}{Example}
\newtheorem{prob}{Problem}\newtheorem{conj}{Conjecture}
\newtheorem{defn}{Definition}
\newtheorem{rem}{Remark}
\newtheorem{ack}{Acknowledgements}
\renewcommand{\theack}{\hskip -4 pt}
\def\cprime{$'$}
\newcommand{\C}{{\mathbb C}} 
\newcommand{\R}{{\mathbb R}} \newcommand{\N}{{\mathbb N}}
\newcommand{\Z}{{\mathbb Z}}
\newcommand{\abs}[1]{{\left| {#1} \right|}} \newcommand{\p}[1]{{\left(
      {#1} \right)}} \newcommand{\jtf}[1]{{#1}^{\diamond}}

\newcommand{\Oh}[1]{{O \p{#1}}}

\renewcommand{\Re}{\operatorname{Re}} \renewcommand{\Im}{\operatorname{Im}}

\author{Johan Andersson\thanks{Department of Mathematics, Stockholm University, 
SE-106 91 Stockholm SWEDEN, Email: {\texttt {johana@math.su.se}}.}}

\date{}
\maketitle

\begin{abstract}We prove a variant of the Mergelyan approximation theorem that allows us to approximate functions that are analytic and nonvanishing in the interior of a compact set $K$ with connected complement, and whose interior is a Jordan domain, with nonvanishing polynomials. This result was proved earlier by the author in the case of a compact set $K$ without interior points, and independently by Gauthier for this case and the case of strictly starlike compact sets. We apply this result on the Voronin universality theorem for compact sets $K$, where the usual condition that the function is nonvanishing on the boundary can be removed. We conjecture that this version of Mergelyan's theorem might be true for a general set $K$ with connected complement and show that this conjecture is equivalent to a corresponding conjecture on  Voronin Universality. 
\end{abstract}

\section{Introduction}

\subsection{Voronin Universality}

Voronin \cite{Voronin2, Voronin0} proved the following Theorem:
\begin{thm}(Voronin)  Let 
 $K=\{ z \in \C :|z-3/4| \leq r\}$ for some $r<1/4$, and suppose that $f$ is any continuous nonvanishing function on $K$ that is analytic  in the interior of $K$. Then 
$$\liminf_{T \to \infty} \frac 1 T \mathop{\rm meas} \left \{t \in [0,T]:\max_{z \in K} \abs{\zeta(z+it)-f(z)}<\varepsilon \right \}>0. $$
\end{thm}

It is well-known that the fact  that $f(z)$ is nonvanishing on $K$ can be relaxed to assuming that $f(z)$ is nonvanishing in the interior of $K$ and may allow zeroes on the boundary, see for example  \cite[p. 251, Theorem 2]{VorKar}. This follows from Theorem 1 in the following way: Consider 

$$f_\xi(z)=f\p{\frac 3 4+(1-\xi)\p{z-\frac 3 4}}, \qquad (0<\xi<1).$$  
If $f(z)$ is nonvanishing in the interior of $K$ then $f_\xi(z)$ is nonvanishing on $K$ for any $\xi>0$. Choose $\xi$ small enough so that $$|f_\xi(z)-f(z)|<\varepsilon/2,$$ for $z \in K$. By Theorem 1 there exist $T$ with positive density so that $$\abs{\zeta(z+iT)-f_\xi(z)}<\varepsilon/2, \qquad (z \in K).$$ The conclusion follows from the triangle inequality. We use a variant of this proof method to prove stronger results later in this paper.

Bagchi \cite{Bagchi}  (see also Steuding \cite[Theorem 1.9]{Steuding}) generalized Theorem 1 to other compact sets than ``little discs''. 
\begin{thm} (Bagchi) Theorem 1 is true when $K$ is any compact set with connected complement lying entirely within $1/2<\Re(s)<1$.
\end{thm}

\subsection{Mergelyan's theorem}

 One important tool needed to prove Theorem 2  is Mergelyan's theorem.
\begin{thm}(Mergelyan) Assume that $K$ is a compact set with connected complement and that $f(z)$ is a function analytic in the interior of $K$ and continuous on $K$. Then there exists for any $\varepsilon>0$ some polynomial $p(z)$ such that
$$
 \max_{z \in K}\abs{f(z)-p(z)}<\varepsilon.
$$
\end{thm}
This was proved in Mergelyan \cite{Mergelyan} and is one of the major theorems in complex approximation theory. For different treatments see Carleson \cite{Carleson} or Rudin \cite[Theorem 20.5]{Rudin}.

\section{Removing nonvanishing on the boundary?}

\subsection{Two conjectures}
 One may ask whether we can still remove the condition that $f(z)$ is nonvanishing on the boundary of $K$ in Theorem 2.  We believe this might be true, but we have not been able to prove this in full generality, so we state this as a conjecture.
\begin{conj} Let $K$ be a compact set with connected complement lying in the strip $1/2<\Re(s)<1$, and  $f(z)$  some continuous function on $K$ that is analytic and nonvanishing in the interior of $K$. Then for any $\varepsilon>0$ we have that
$$\liminf_{T \to \infty} \frac 1 T \mathop{\rm meas} \left \{t \in [0,T]:\max_{z \in K} \abs{\zeta(z+it)-f(z)}<\varepsilon \right \}>0. $$
\end{conj}
Conjecture 1 is related to the following conjectured variant of Mergelyan's theorem.
\begin{conj}  
  Assume that $K$ is a compact set with connected complement and that $f(z)$ is a continuous function on $K$  that is analytic and nonvanishing in the interior of $K$. Then there exists for any $\varepsilon>0$ some polynomial $p(z)$ that is nonvanishing on $K$ such that
\begin{gather*}
 \max_{z \in K}\abs{f(z)-p(z)}<\varepsilon.
\end{gather*}
\end{conj}
\begin{rem} Gauthier informed the author that he had thought about this problem in the seventies, although the problem itself is not published. It is related to  results in Gauthier-Roth-Walsh \cite{GauthierRothWalsh}.
\end{rem}

\subsection{Relating Mergelyan's theorem and Voronin universality}

\begin{thm} Conjecture 1 and Conjecture 2 are equivalent. \end{thm}
\begin{proof} $i)$ {\em Conjecture 2 implies Conjecture 1.} We employ the same argument as in the proof of Theorem 2 in \cite{Andersson}.  By Conjecture 2 we can approximate $f(z)$ by a polynomial $p(z)$ such that \begin{gather} \label{ui} \abs{p(z)-f(z)}< \varepsilon/2,  \qquad \qquad (z \in K), \end{gather}
where $p(z)$ is nonvanishing on $K$. By Theorem 2  we  have that  
$$\liminf_{T \to \infty} \frac 1 T \mathop{\rm meas} \left \{t \in [0,T]:\max_{z \in K} \abs{\zeta(z+it)-p(z)}<\varepsilon/2 \right \}=\delta>0. $$
From the inequality
$$\max_{z \in K} \abs{\zeta(z+it)-p(z)}<\varepsilon/2,$$
it follows by  the  triangle inequality and \eqref{ui} that 
$$\max_{z \in K} \abs{\zeta(z+it)-f(z)}<\varepsilon.$$
Hence 
$$\liminf_{T \to \infty} \frac 1 T \mathop{\rm meas} \left \{t \in [0,T]:\max_{z \in K} \abs{\zeta(z+it)-f(z)}<\varepsilon \right \}\geq\delta>0. $$ \qed

 $ii)$ {\em Conjecture 1 implies Conjecture 2.}  
   Since $K$ is compact we can choose $\varepsilon$ sufficiently small so  we have that 
   \begin{gather*} K_0=3/4+\varepsilon K  \subset \{z:|z-3/4| \leq 1/8 \}. \\ \intertext{Let} 
g(z)=f\p{\frac{z-3/4} \varepsilon}.
\end{gather*}
It is clear that $g(z)$  is analytic and nonvanishing in the interior of $K_0$ by the fact that $f(z)$  is analytic and nonvanishing in the interior of $K$.   Since $K_0$ lies strictly in $5/8\leq \Re(s)\leq 7/8$ we have that 
$$\sigma=\inf_{z \in K_0} \Re(z) \geq 5/8.$$ 
Standard zero-density estimates for the Riemann zeta-function, for example the estimate of Ingham \cite{Ingham} (see also Ivi{\'c} \cite[Chapter 11]{Ivic})
\begin{gather*}
  N(\sigma,T)\ll  T^{3(1-\sigma)/(2-\sigma)}\log^5 T,
 \end{gather*}
 where $$N(\sigma,T)=\sharp{} \{s:\zeta(s)=0, \Re(s) \geq \sigma, |\Im(s)| \leq T  \},$$  
denote the number of zeroes\footnote{The Riemann hypothesis says that $N(\sigma,T)=0$ for $\sigma>1/2$} of the Riemann zeta function in a rectangle  implies that $N(5/8,T)\ll_\varepsilon T^{9/11+\varepsilon}$ for any $\varepsilon>0$ and in particular that
$$\lim_{T \to \infty} \frac 1 T \mathop{\rm meas} \left \{t \in [0,T]:\exists z \in K_0:\zeta(z+it)=0\right \}=0. $$
In contrast, we have positive density in Conjecture 1. This means that for any $\varepsilon>0$ we can find some $T \geq 2$ such that the Riemann zeta-function $\zeta(z+iT)$ has no zeroes on $K_0$, i.e.
\begin{gather} \label{Eq0}
  \min_{z \in K_0} \abs{\zeta(z+iT)} = \delta>0, \\ \intertext{and that}
 \label{Eq1}
  \max_{z \in K_0}\abs{g(z)-\zeta(z+iT)}<\frac \varepsilon 2.
\end{gather}
We have that $\zeta(z+iT)$ is an analytic function for $|z| \leq 1$   and it can thus be approximated by a polynomial $q(z)$ such that
\begin{gather} \label{Eq2}
  \abs{ q(z) -\zeta(z+iT)} <\min (\varepsilon/2,\delta/2), \qquad \qquad (\abs{z} \leq 1).
\end{gather}
In particular this is true for $z \in K_0$ and by combining equations \eqref{Eq1} and \eqref{Eq2} we find that
\begin{gather} \label{Eq88}
  \abs{g(z)-q(z)}<\varepsilon, \qquad (z \in K_0).
\end{gather}
where $q(z)$ is a polynomial that by \eqref{Eq0} and \eqref{Eq2} fulfils
 \begin{gather*} 
  \abs{q(z)}\geq \frac \delta 2 >0, \qquad (z \in K_0),
\end{gather*}
and is thus nonvanishing on $K_0$. Let
\begin{gather*}
  p(z)=q\p{3/4+\varepsilon z}.
\end{gather*}
By the construction of the set $K_0$, the function $g(z)$,  Eq. \eqref{Eq88} and the fact that $q(z)$ is nonvanishing on $K_0$   it is clear that the polynomial $p(z)$ is nonvanishing on $K$ and that
\begin{gather*}
  \sup_{z \in K} \abs{f(z)-p(z)}<\varepsilon.
\end{gather*} 
 \end{proof}
\begin{rem}
  Universality theorems are known for many different Dirichlet series, including the Selberg class, see Steuding \cite{Steuding}. Conjecture 1 can be formulated for any element of this class as well and we still have equivalence between Conjecture 1 and Conjecture 2 by the same proof method. Kaczorowski-Perelli \cite{KacPer} have proven a suitable zero-density-estimate replacing Ingham's and together with Steuding's universality results the same proof holds. The only complication is that we might need to move the set $K_0$ closer to the line $\Re(s)=1$ than being centered at $3/4$, since the universality results holds in a more narrow strip. One interesting consequence of this is that a universality result of the same type as Conjecture 1 (possible for a more narrow strip) for one function in the Selberg class implies the same type of result for any other element in the Selberg class.
\end{rem}

\begin{rem}
  In case we know the Riemann hypothesis for an $L$-function we do not need positive density in Conjecture 1 and we still have equivalence between Conjectures 1 and 2. While we do not know the Riemann hypothesis for any element in the Selberg class, recent important  results of Drungilas-Garunk{\v{s}}tis-Ka{\v{c}}enas \cite{DruGarKac} proves the Voronin universality theorem in the strip $0.848\ldots <\Re(s)<1$ of the  Selberg zeta-function for the full modular group. In this case the Riemann hypothesis is known to hold. We can then formulate a weaker version of Conjecture 1 and still prove that it implies Conjecture 2.
\end{rem}

\section{Proof of our conjectures for special cases}

In \cite{Andersson} we managed to show Conjecture 1 and 2 for the case of compact sets without interior points\footnote{This special case of Mergelyan's theorem is called  Lavrent{\cprime}ev's theorem. For a different proof of this result see Gauthier \cite[Proposition 32]{Gauthier}.}. When applied on the Voronin theorem it simplifies the statement, since not only the assumption that $f(z)$ is nonvanishing on $K$ can be removed completely, but also the assumption that $f(z)$ is analytic on the interior of $K$ can be removed since the interior of $K$ is empty. This allowed us to prove a criterion of Bagchi in this special case. In contrast, even if we manage to prove Conjecture 2 for a general compact set $K$ it will not imply anything similar. This is because while the condition that $f(z)$ is nonzero on the boundary of $K$ might be removed (if conjectures 1 and 2 are true) it is easy to see that the condition that $f(z)$ is nonzero in the interior of $K$ cannot be removed.

While we can not treat the general case of Conjectures 1 and 2, we have managed to show some partial results.
\begin{thm}
  Conjecture 1 and 2 are true if the interior of $K$ is a Jordan domain.
\end{thm}
We remark that a Jordan domain is an open connected set that is bounded by a Jordan curve, see e.g. Palka \cite[p. 34]{Palka}.

\begin{proof}
It is sufficient to prove Conjecture 2 for these compact sets $K$ since Conjecture 1 will follow from Conjecture 2 in this special case, in the same way as in the general case, see the first part of the proof of Theorem 4.

 Let  $O=K^o$ be a Jordan domain.  By the Carath{\'e}odory-Osgood-Taylor theorem\footnote{Problem suggested by Osgood and proved independently by  Carath{\'e}odory  \cite{Cara} and Osgood-Taylor \cite{Osgood}. For text book references, see  Palka \cite[Theorem 4.9]{Palka} or Rudin \cite[Theorem 14.19]{Rudin}.}, the  Riemann mapping theorem  that there  exists a holomorphic bijection $\phi:D\to O$ between the disc $ D=\{z \in \C: |z|<1\}$  and $O$, can be extended to a continuous map $\phi:\overline D \to \overline{O}$.
It is clear that
$$f(z)=f(\phi(\phi^{-1}(z)) $$
on $\overline{O}$. Choose
$$H(z)=f(\phi((1-\xi)\phi^{-1}(z))), $$
for a sufficiently small $\xi$ such that
\begin{gather} \label{uttt} |H(z)-f(z)| \leq \varepsilon/3, \qquad \qquad (z \in \overline{O}). \end{gather} 
Tietze's extension theorem \cite[Theorem 20.4]{Rudin} allows $H$ to be extended  to a continuous function on $K$ such that
\begin{gather} \label{ut2}
\abs{H(z)-f(z)} \leq \varepsilon/3, \qquad \qquad (z \in K).
\end{gather} 
 By the construction it is clear that $H(z)$ is continuous on $K$, nonvanishing on $\overline{O}$ and analytic on $O$. Thus
\begin{gather*}
  \sup_{z \in \overline{O}} \abs{H(z)}=\delta>0,
\end{gather*} 
 and by Theorem 3 we can choose a polynomial $P(z)$ such that
\begin{gather} \label{ut3}
|P(z)-H(z)| < (\varepsilon/3,\delta/2), \qquad \qquad (z \in K).
\end{gather}
By the triangle inequality it is clear that 
$$|P(z)|> \delta/2, \qquad \qquad (z \in \overline{O}).$$ 
We now use the same proof method as in the proof of Theorem 1 in \cite{Andersson}. Let  
\begin{gather}
\notag P(z)=c_0 \prod_{k=1}^m (z-z_k), \\ \intertext{where $z_k$ denote the zeros of $P(z)$. Since the polynomial $P(z)$ has no zeros in the interior of $K$, each zero $z_k$ must lie in the boundary of $K$ or outside $K$ and there exist  sequences $z_{k,n}$ of points in $\C \setminus K$ such that $\lim_{n \to \infty} z_{k,n}=z_k$. Define}
 \notag  p_n(z)=c_0 \prod_{k=1}^m (z-z_{k,n}). \\ \intertext{
Since  all the coefficients of the polynomials $p_n(z)$ will converge to the coefficients of $P(z)$, it is clear that $p_n(z)$ will converge to $P(z)$ uniformly on the compact set $K$. Hence there exists an $n$ such that}
\label{u2} |p_n(z)-P(z)|<\varepsilon/3.
\end{gather}
 Since  $z_{k,n}$ denote points in $\C \setminus K$, the polynomial $p_n(z)$ will have all its zeroes outside of $K$ and the polynomial will be nonvanishing on $K$. We can therefore choose $p(z)=p_n(z)$. By the triangle inequality and the inequalities  \eqref{ut2}, \eqref{ut3} and  \eqref{u2} we obtain  the inequality $|p(z)-f(z)|<\varepsilon$ for every $z \in K$.
\end{proof}

We will also be able to treat the case of finitely many interior components in some special cases.
\begin{thm}
  Conjectures 1 and 2 are true if $K$ has finitely many maximally connected open subsets $O$, each such subset is a Jordan domain, and
furthermore if $O_1 \neq O_2$ are two such subsets then $f(z)$ is nonvanishing on $\overline{O_1} \cap \overline {O_2}$.
\end{thm}

\begin{proof}
  As in the proof of Theorem 5 it is sufficient to prove Conjecture 2 for these sets, since Conjecture 1 is a consequence. 
Also, it is sufficient to construct a continuous function $H(z)$ that is analytic in the interior of $K$, which is denoted by $K^o$ and nonzero on its closure $\overline{K^o}$, such that
\begin{gather} \label{aj}
 \abs{f(z)-H(z)}<\varepsilon/3, \qquad  \qquad (z \in \overline{K^o}),
\end{gather}
since the rest of the proof follows in the same way as from  equation \eqref{uttt} in the proof of Theorem 5.  
Since $K$ has finitely many maximal connected open sets, we have  the decomposition
$$
 \overline{K^o}=\bigcup_{l=1}^L K_l,
$$
as a disjoint union of connected compact sets $K_l$. For each such set $M= K_l$ we have that
\begin{gather} \label{ajaj} M^o=\bigcup_{n=1}^{N} O_{n}, \end{gather}
where $O_n$ are disjoint Jordan domains. 
Now form the graph with vertices $\{O_n:n=1,\ldots, N\}$ and an edge between $n$ and $k$ if and only if $\overline{O_n} \cap \overline{O_k} \neq \emptyset$ and $n \neq k$. It is clear that the graph is connected from the fact that $M$ is connected. Furthermore it is clear that the graph is a tree, since if we have a cycle in the graph we can construct
a Jordan curve $J \subset M \subset K$ such that there exist points both on the inside and the outside of the curve that are not in the set $K$.  This Jordan curve can be chosen as the union of {\em half}\footnote{{\em half} here means that the two intersection points (one intersection point with the closure of each neighbour in the cycle) on the Jordan curve that bounds the Jordan domain will divide it into two curves and we choose either one of them.} of the boundary for each Jordan domain in the cycle. 
By the Jordan curve theorem, this  violates our assumption that the complement of $K$ is connected.

 Since  the graph is a tree, we can choose some root $O_{\alpha_1}$ of the tree to obtain a rooted tree. Any rooted tree induces a partial ordering, the tree-order where $u \leq  v$ if  and only if the unique path from the root to $v$ passes through $u$. Furthermore by topological sorting any partial order admits a total order
\begin{gather*}
         O_{\alpha_1}< O_{\alpha_2}< \cdots <O_{\alpha_N}.
\end{gather*}
With this ordering we will always have 
\begin{gather*}
\overline{O_{\alpha_n}} \cap\overline{O_{\alpha_k}} \neq \emptyset,
\end{gather*}
 for exactly one $k<n$,  in fact the vertex $O_{\alpha_k}$ in the rooted tree will be the parent of the vertex  $O_{\alpha_n}$.  Again, $\overline{O_{\alpha_n}} \cap \overline{O_{\alpha_k}}$ cannot contain more than one point, since then we can construct a Jordan curve violating the fact that we know that  the complement of $K$ is connected. Thus the intersection consists of exactly one point.
Let 
\begin{gather*}
\{z_n\}=\overline{O_{\alpha_n}}  \cap \overline{O_{\alpha_k}},
\end{gather*}
for  this $k$. We see that
\begin{gather} \label{p}
  \delta=\min_{2 \leq n \leq N} \abs{f(z_n)} = \min_{z \in \overline{O_{\alpha_n}} \cap \overline{O_{\alpha_k}}, n \neq k} \abs{f(z)}>0.
\end{gather}
Also let
\begin{gather} \label{pp}
  C=\max_{z \in K} \abs{f(z)}+1.
\end{gather}
By Theorem 5 we now choose for  each $O_{\alpha_n}$ a function $f_n$ so that
\begin{gather} \label{ppp}
  \max_{z \in \overline{O_{\alpha_n}}}\abs{f_n(z)-f(z)}< \xi = \min\left(1,\delta/2,  \p{\delta/(3C)}^N \varepsilon/3\right),
\end{gather}
where $f_n(z)$ is analytic on $O_{\alpha_n}$, and nonzero and continuous on $\overline{O_{\alpha_n}}$. We will now  glue the functions together.
Define recursively
\begin{gather} \label{Hdef} 
 \begin{split}
  H_1(z)&= f_1(z), \qquad  \qquad  \hskip 1pt \, \, \,  \, \, \, \text{for}   \qquad z \in \overline{O_{\alpha_1}},\\
    H_n(z)&= \begin{cases} H_k(z), & \text{for} \qquad   z \in \overline{O_{\alpha_k}}, \qquad   1 \leq k \leq n-1, \\
    \frac{H_{n-1}(z_n) f_n(z)}{f_{n}(z_n)},  & \text{for}   \qquad z \in \overline{O_{\alpha_n}}, \qquad 2 \leq n \leq N.
 \end{cases}
\end{split}
\end{gather}
It is clear that $H(z)=H_N(z)$ is continuous on $M$ since we have ensured continuity at the points $z_n$. By \eqref{Hdef} we have that
\begin{gather*}
  H_n(z)-f(z)
  =\frac{((H_{n-1}(z_n)-f(z_n))+(f(z_n)-f_n(z_n)))f_n(z)}{f_n(z_n)} + f_{n}(z)-f(z), \\ 
\intertext{when $z \in \overline{O_{\alpha_n}}$. By \eqref{p}, \eqref{pp}, \eqref{ppp} and the triangle inequality it is clear that $|f_n(z)|<C$ and $|f_n(z_n)|>\delta/2$, and if we define} 
\xi_n=  \max_{z \in \overline{O_{n}}} \abs{H_n(z)-f(z)},  \\ \intertext{another application of the triangle inequality and \eqref{ppp} gives us that}
  \xi_n <  \left( \max_{1\leq k \leq n-1}\xi_{k}+\xi \right) \frac{2C} {\delta} +\xi, \qquad (2 \leq n \leq N). 
\end{gather*}
By \eqref{ppp} and \eqref{Hdef} we have $\xi_1 < \xi$ and by a simple induction argument using the fact that $C/\delta \geq 1$ it follows that $\xi_n < (3C/\delta)^n \xi$ for $1 \leq n \leq N$. By the choice of $\xi$ in \eqref{ppp} we see that $\xi_n < \varepsilon/3$ for $1 \leq n \leq N$ and by defining $H(z)$ on each $M=K_l$ for $1 \leq l \leq L$ as above we obtain the inequality \eqref{aj}.  \end{proof}

\begin{rem}
 I started to think about this problem in August 2009, and most of the results in this paper are from that month. I did initially choose to publish just the empty interior case \cite{Andersson}, since it allowed an especially nice formulation of the Voronin universality theorem. The decision to proceed with the publication of these results  was partly inspired by seeing a copy of the paper  \cite{Gauthier} of Paul M. Gauthier, where he treats the strictly starlike and empty interior case (independently from the present author). 
\end{rem}

\section{Open problems and further research}

While we believe Conjectures 1 and 2 {\em might} be true, we do not have very strong reasons to believe in it. In fact, some quite strange sets can be constructed. If we consider a simply connected open set that is not bounded by a Jordan curve, the situation is more complicated, see the literature on the Carath{\'e}odory's theory of prime ends (see e.g. \cite{Cara} or \cite{Epstein}).

A nontrivial example suggested by Anthony G. O'Farrell is the Cornucopia set\footnote{This set can be found in Gamelin \cite{Gamelin} and has also found other applications in the theory of complex approximation by polynomials, see O'Farrell and  P{\'e}rez-Gonz{\'a}lez \cite{Farrell}.}: Let us have a compact set with two interior components. One open disc and one strip. Let the strip go around the disc indefinitely and successively approach the disc at the same time as it thins out.  Then the strip will not be a Jordan domain, because we have that the boundary of the disc is in fact a subset of the boundary of the strip. Since interior points in the disc will not belong to the closure of the strip, this means that the closure of the strip will not be simply connected.  In contrast the closure of any Jordan domain is simply connected.

Also the case when we have infinitely many maximal connected interior sets seems difficult, even when all the sets are Jordan domains.  Difficult configurations can be found, when they touch each other, are close to each other and when they look like snakes, i.e. even if the area tends to zero, their length stays the same.  

Of course if the conjectures are false, it would be interesting to have a counterexample.  In any case it would be interesting to have more cases where the conjectures are known to be true. Some cases seems easier. Certain cases with an infinite number of open components can certainly be considered. Also it seems likely that the condition that $f(z)$ is nonzero on the intersection between closures of disjoint maximal connected interior sets can be removed.

\begin{ack}
 The author is grateful to Lawrence Zalcman, Anthony G. O'Farrell, Maria Roginskaya and Paul M. Gauthier for expressing interest in this problem, and giving me motivation to finish this paper. The author is also grateful to the referee for pointing out  problems with previous versions of the paper.
\end{ack}

\bibliographystyle{plain}

\end{document}